\documentclass[conference,letterpaper]{IEEEtran}

\usepackage[utf8]{inputenc} 
\usepackage[T1]{fontenc}
\usepackage{url}
\usepackage{ifthen}
\usepackage{cite}
\usepackage[cmex10]{amsmath} 

\usepackage{epsfig}
\usepackage{times}
\usepackage{float}
\usepackage{afterpage}
\usepackage{amsmath}
\usepackage{amstext}
\usepackage{amssymb,bm,bbm}
\usepackage{latexsym}
\usepackage{color}
\usepackage{graphicx}
\usepackage{amsmath}
\usepackage{amsthm}
\usepackage{graphicx}
\usepackage[center]{caption}
\usepackage{subfigure}
\usepackage{booktabs}
\usepackage{multicol}
\usepackage{lipsum}
\usepackage{dblfloatfix}
\usepackage{mathrsfs}
\usepackage{cite}
\usepackage{tikz}
\usepackage{pgfplots}
\pgfplotsset{compat=newest}
 \usepackage{enumitem}
\usepackage{makecell}

\allowdisplaybreaks

\newcommand{\rmd}{\mathrm{d}}
\newcommand{\bbE}{\mathbb{E}}\newcommand{\rme}{\mathrm{e}}

\newcommand{\bbN}{\mathbb{N}}

\newcommand{\bbR}{\mathbb{R}}

\newcommand{\bfA}{\mathbf{A}}

\newcommand{\bfp}{\mathbf{p}}

\newcommand{\cW}{\mathcal{W}}

\DeclareMathOperator*{\argmin}{argmin}

\theoremstyle{mystyle}
\newtheorem{theorem}{Theorem}
\theoremstyle{mystyle}
\newtheorem{lemma}{Lemma}
\theoremstyle{mystyle}
\theoremstyle{mystyle}
\newtheorem{corollary}{Corollary}
\theoremstyle{mystyle}
\theoremstyle{remark}
\theoremstyle{mystyle}
\theoremstyle{mystyle}
\theoremstyle{mystyle}
\theoremstyle{discussion}
\theoremstyle{mystyle}
\theoremstyle{mystyle}

 \usepackage{enumitem}

\IEEEoverridecommandlockouts
\begin{document}

\title{
       Linearity-Inducing Priors for Poisson Parameter Estimation Under \(L^{1}\) Loss   \thanks{
This work was supported in part by the U.S. National Science Foundation under Grant ECCS-2335876.} }
 \author{%
  \IEEEauthorblockN{Leighton P. Barnes\IEEEauthorrefmark{1},
                    Alex Dytso\IEEEauthorrefmark{2}, 
                    and H. Vincent Poor\IEEEauthorrefmark{3}}
  \IEEEauthorblockA{\IEEEauthorrefmark{1}%
                    Center for Communications Research,
                    Princeton, NJ 08540, USA,
                    l.barnes@idaccr.org}
   \IEEEauthorblockA{\IEEEauthorrefmark{2}%
                    Qualcomm Flarion Technology, Inc.,
                    Bridgewater, NJ 08807, USA,
                    odytso2@gmail.com}
  \IEEEauthorblockA{\IEEEauthorrefmark{3}%
                    Princeton University, 
                    Princeton, NJ 08544,  USA,
                    poor@princeton.edu}

}                  

\maketitle

\begin{abstract}
We study prior distributions for Poisson parameter estimation under $L^1$ loss. Specifically, we construct a new family of prior distributions whose optimal Bayesian estimators (the conditional medians) can be any prescribed increasing function that satisfies certain regularity conditions. In the case of affine estimators,  this family is distinct from the usual conjugate priors, which are gamma distributions. Our prior distributions are constructed through a limiting process that matches certain moment conditions. These results provide the first explicit description of a family of distributions, beyond the conjugate priors, that satisfy the affine conditional median property; and more broadly for the Poisson noise model they can give any arbitrarily prescribed conditional median.
\end{abstract}

\section{Introduction}

In Bayesian estimation theory, the optimal estimator under $L^1$ loss is given by the conditional median. In particular, suppose that a random variable $X$ is drawn from a prior distribution $P_X$, and we are given a noisy measurement $Y$ of $X$ that is drawn from some noise model $P_{Y|X}$. We use $\widehat{X}(Y)$ to denote our estimate of $X$ from $Y$. The conditional median is optimal in the following sense:\footnote{ In case there are multiple optimizers, which can be the case for discrete distributions, the median is always taken to be $\mathsf{med}(X|Y=y) = F_{X|Y=y}^{-1}(\frac{1}{2})$ where $F_{X|Y=y}^{-1}(p) \in (0,1)$ is the left-inverse of the conditional cumulative distribution function (cdf) also known as the quantile function.}
\begin{equation} \mathsf{med}(X|Y=\cdot) = \argmin_{\widehat{X}} \mathbb{E} \left[|\widehat{X}(Y) - X| \right],
\end{equation}
where the argmin is over all possible measurable functions and the expectation is jointly over $X,Y$.

Our focus is on the Poisson noise model, where the observation $Y$ is a non-negative integer with 
\begin{equation}
     P_{Y|X}(y|x) = \frac{1}{y!} x^y \rme^{- x} , x \ge 0,\, y\in \bbN_0, \label{eq:Poisson_pmf}
 \end{equation}
where $\bbN_0$ denotes the set of non-negative integers, and we use the convention  that  $0^0=1$. The Poisson parameter $X$ is supported on the non-negative reals. 

The Poisson noise model is fundamental to applications involving counting discrete events, such as the number of photons in optical communications \cite{mceliece1979practical,shamai1990capacity, verdu1999poisson,lapidoth2008capacity,dytso2021properties} or particles in molecular communications \cite{Farsad2020}, and the number of neural spikes in neuroscience \cite{neuro1, Shadlen3870}. The latter application especially establishes Poisson noise, and understanding its differences from other neural noise models such as Gaussian noise, is a topic of great interest in understanding noisy large-scale neural systems and learning models.

For $L^2$ error, where the optimal estimator is the conditional mean instead of the conditional median, and it is known that
\begin{equation}\label{eq:linearconditionalmean}
    \mathbb{E}[X|Y=y] = ay+b \; ,
\end{equation}
i.e., the optimal estimator is an affine function, if and only if  $0<a <1, b>0$  and the prior distribution $P_X$ is a gamma distribution \cite{dytso2020estimation} with the probability density function
\begin{equation}
    f_X(x) = \frac{\alpha^{\theta}}{\Gamma \left(  \theta \right) } x^{\theta-1} \rme^{-\alpha x}, \, x \ge 0, 
\end{equation}
where 
\begin{equation}
    \alpha =  \frac{1-a}{a} , \, \theta  =  \frac{b}{a} .  \label{eq:parameters}
\end{equation}
The gamma distribution is the conjugate prior for this exponential family; and for a natural exponential family, the optimal $L^2$ estimator for the so-called \emph{mean parameter} (which is just $X$ in the Poisson case) is affine if and only if the prior distribution on $X$ is the conjugate prior \cite{diaconis1979conjugate,chou2001characterization}.

In this work, we consider the corresponding property for the conditional median, and ask if there are prior distributions $P_X$ such that
\begin{equation} \mathsf{med}(X|Y=y) = ay + b  \; .\label{eq:affineconditionalmedian}\end{equation}
This follows a corresponding line of research that the authors have undertaken for the additive Gaussian noise model and some other exponential families \cite{BDP_ISIT2023,BDLP_IT2024,BDLP_ISIT2024}; see also \cite{akyol2012conditions} for similar results. In the Gaussian case, there is a unique prior that guarantees a linear conditional median, and it corresponds to the conjugate prior (which is another Gaussian with a particular variance).

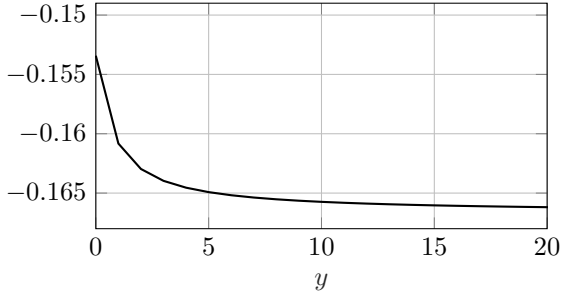
\begin{figure}
    \centering  
%
%
\definecolor{mycolor1}{rgb}{0.00000,0.44700,0.74100}%
\begin{tikzpicture}

\begin{axis}[%
width=6cm,
height=3cm,
at={(1.011in,0.642in)},
scale only axis,
xmin=0,
xmax=20,
xlabel style={font=\color{white!15!black}},
xlabel={$y$},
ymin=-0.168,
ymax=-0.149,
yticklabel style={
            /pgf/number format/.cd,
            fixed,
            precision=5
        },
axis background/.style={fill=white},
xmajorgrids,
ymajorgrids,
legend style={legend cell align=left, align=left, draw=white!15!black}
]
\addplot [color=black, thick]
  table[row sep=crcr]{%
0	-0.153426409720027\\
1	-0.16082650499167\\
2	-0.16296984313822\\
3	-0.163969625574552\\
4	-0.164545558602008\\
5	-0.164919405643965\\
6	-0.165181462725114\\
7	-0.165375278749598\\
8	-0.165524407814813\\
9	-0.165642692642934\\
10	-0.165738798081843\\
11	-0.165818423477616\\
12	-0.165885470630682\\
13	-0.165942700327551\\
14	-0.165992120834603\\
15	-0.166035227841375\\
16	-0.166073158130191\\
17	-0.166106791105769\\
18	-0.166136817714158\\
19	-0.166163788347164\\
20	-0.16618814684475\\
21	-0.166210255084026\\
22	-0.166230411060683\\
23	-0.166248862386841\\
24	-0.16626581650579\\
25	-0.166281448518792\\
26	-0.166295907250744\\
27	-0.166309320000241\\
28	-0.166321796295307\\
29	-0.166333430889384\\
30	-0.166344306171041\\
31	-0.166354494117009\\
32	-0.166364057886394\\
33	-0.166373053130716\\
34	-0.166381529077132\\
35	-0.166389529429395\\
36	-0.166397093121322\\
37	-0.16640425495024\\
38	-0.166411046112071\\
39	-0.166417494655448\\
40	-0.166423625868884\\
41	-0.166429462612054\\
42	-0.166435025600578\\
43	-0.1664403336516\\
44	-0.166445403896397\\
45	-0.166450251964935\\
46	-0.166454892147158\\
47	-0.166459337533055\\
48	-0.166463600136353\\
49	-0.166467691002818\\
50	-0.166471620305813\\
};

\end{axis}
\end{tikzpicture}%
    \caption{Plot of the conditional median in expression \eqref{eq:Conditional_median_poisson}  minus the conditional mean in \eqref{eq:linearconditionalmean} for $a=b=\frac{1}{2}$. The difference decays like $O(1/y)$ as shown in \cite{chen1986bounds}.}
    \label{fig:comp_median_to_mean}
\end{figure}

In contrast to both the $L^2$ error case and the Gaussian noise model case, the gamma distribution does not induce an affine conditional median under Poisson noise since
 \begin{equation}
\mathsf{med}(X|Y=y)= \frac{1}{\alpha} \gamma^{-1} \left(\frac{1}{2}, \theta + y \right), \label{eq:Conditional_median_poisson}
 \end{equation}
 where $\gamma^{-1}$ is the inverse of the lower incomplete gamma function and needs to be computed numerically. 
 However, it is a near miss, as shown in Fig.~\ref{fig:comp_median_to_mean}; see \cite{chen1986bounds} for a precise error analysis where first- and second-order error terms are derived. 

\subsection{Paper Outline and Contribution}
Our first  contribution is to show that for  any increasing function $f$ (that satisfies a particular regularity condition) there exists a prior such that
\begin{equation}
\mathsf{med}(X|Y=y) = f(y), \, y \in \bbN_0.
\end{equation}
Specializing this result to affine functions $f$, we demonstrate a new family of prior distributions that \emph{do} induce an affine conditional median. They are distinct from the usual conjugate prior/gamma distribution, but may have a qualitatively similar tail as seen in Fig.~\ref{fig:cdfs}. They are constructed by a limiting procedure of matching particular moment conditions, and this limiting procedure is guaranteed to have at least a subsequence of distributions that converge in total variation to a measure with the desired properties.

To the best of our knowledge, this is the first explicit description of a distribution that is distinct from the usual conjugate prior, but satisfies the affine conditional median property \eqref{eq:affineconditionalmedian}. In this sense, it could be viewed as a family of \emph{$L^1$-conjugate} priors, but without necessarily preserving the property that conditioning is closed with respect to the family of distributions. 

We proceed as follows: In Section \ref{sec:main} below we state our main existence results, and then Section \ref{sec:proofs} is devoted to proving these results by constructing our distributions via a limiting process. In Section \ref{sec:Examples} we conclude by showing examples and discussing the implications of this construction.

\section{Main Results} \label{sec:main}
In this section, we present our main results. We begin by presenting the existence result, followed by a description of a procedure that can be used to construct such a distribution. 

\begin{theorem} \label{thm:main_theorem}  Suppose that $f:\bbN_0 \to \bbR_{+}$  satisfies the following properties: 
\begin{itemize}
    \item  $f$ is an increasing function up to a point  $c_0$ and constant afterwards (i.e., $f(i-1) < f(i), \, i \le c_0$ and $f(i+1)=f(c_0), i \ge 
 c_0$). The point $c_0$ can be infinity; and
    \item there exists $\kappa \ge 1$ such that 
    \begin{equation} \label{eq:sum_cond}
\sum_{i=1}^{c_0} \left(  \frac{f \left( \lfloor \frac{i}{\kappa}  \rfloor  \right)  }{f(i)} \right)^{ \frac{i}{\kappa}    } \rme^{f(i)}  <\infty. 
\end{equation}
\end{itemize}
Then, there exists a distribution $P_X$ supported on 
\begin{equation} \label{eq:arb_support}
    \cW =\{w_i: w_i=f(i-1), \, i \in \bbN \}
\end{equation}  
such that for all $y \in \bbN_0$ \begin{equation}
\mathsf{med}(X|Y=y) = f(y). 
\end{equation}
\end{theorem}

The first assumption in Theorem~\ref{thm:main_theorem} is not very restrictive in view of the following facts: 
\begin{itemize}[leftmargin=*]
\item The assumption that $f$ is increasing is justified by the fact that for  Poisson noise  optimal Bayesian estimators are non-decreasing \cite{nowak2013monotonicity}.  
\item If $X \le A$ has bounded support, then $\mathsf{med}(X|Y=y) \le A, y\in \bbN$, i.e., the conditional median can increase only up to a point. 
\end{itemize}

Specializing, the above result to linear estimators we have the following result. 

\begin{corollary} \label{cor:affine}
   Given $0 < a < \frac{1}{\rme} \approx 0.3678$ and $b>0$,  let  
\begin{equation}
 \cW = \left\{w_i:   w_i = a (i-1) + b, \, i \in \bbN  \right\}.  \label{eq:support}
\end{equation}
Then,  there exists $P_X$ supported on $\cW$  such that for all $y \in \bbN_0$ \begin{equation}
\mathsf{med}(X|Y=y) = ay+b. 
\end{equation}
\end{corollary}
\begin{IEEEproof}
    Setting $f(y) = ay+b$, the sum in \eqref{eq:sum_cond} reduces to 
    \begin{equation}
\inf_{\kappa \ge 1} \sum_{i=1}^{\infty} \left(  \frac{a\lfloor \frac{i}{\kappa}  \rfloor  +b}{a i+b} \right)^{ \frac{i}{\kappa}    } \rme^{a i+b}, 
\end{equation}
which converges provided that $a<\frac{1}{\rme}$. 
\end{IEEEproof}

Note that unlike for the conditional mean case, the distribution inducing a linear conditional median is not necessarily unique. It is also worth pointing out that the range of $a \in \left(0, \frac{1}{\rme}\right)$ in Corollary \ref{cor:affine} is most likely not exhaustive. The full set of admissible values of $a$ is left as an open problem.

\section{Proof } \label{sec:proofs}

\subsection{Equivalent Moment Condition}
We first produce an equivalent moment condition, which is easier to work with.  
\begin{lemma}\label{lem:Moment_Version} Given $f: \bbN_0 \to \bbR_{+}$,  suppose there exists a distribution $P_W$ such that 
\begin{equation}
    \bbE[W^y 1_{W \le f(y)}] = \bbE[W^y 1_{W > f(y)}], \, y  \in \bbN_0. \label{eq:moment_condition}
    \end{equation}
Then, 
\begin{equation}
\rmd P_X(w)= \rmd P_W(w) \rme^{w}
\end{equation}
induces
\begin{equation}
\mathsf{med}(X|Y=y) = f(y), \, y \in \bbN_0,
\end{equation}
provided that $P_X$ is a proper distribution. 
\end{lemma}
\begin{IEEEproof}
An equivalent condition for the prescribed conditional median is
given by the conditional cdf: for all $y \in \bbN_0$
\begin{equation}
    \bbE \left[ 1_{ X \le f(Y)} | Y =y \right]   =\frac{1}{2}
\end{equation}
which is equivalent to:   for all $y \in \bbN_0$
\begin{equation}
    \bbE [ 1_{ X \le f(Y)} | Y =y]   =  \bbE [ 1_{ X > f(Y)} | Y =y] . \label{eq:mid_step}
\end{equation}
Now by defining 
\begin{equation}
\rmd P_W(w) \propto \rme^{-w} \rmd P_X(w),
\end{equation} 
and using the structure of the $P_{Y|X}$ in \eqref{eq:Poisson_pmf} the expression in \eqref{eq:mid_step} can be re-written as
\begin{equation}
\int  w^y  1_{ w  \le f(y) } \rmd P_W(w)  =   \int w^y  1_{ w  > f(y) }  \rmd P_W(w) . 
\end{equation} 
This concludes the proof. 
\end{IEEEproof}

We now proceed to show  the following theorem.

\begin{theorem} \label{thm:moment_theorem} Suppose that  $f: \bbN_0 \to \bbR_0$ is  increasing up to a point $c_0$, such that for every $y \in \bbN_0$
\begin{equation}
   \inf_{ \kappa \ge 1}  \lim_{N \to \infty}  \sum_{i= N+1}^{c_0}   \left(  \frac{f \left( \lfloor \frac{i}{\kappa}  \rfloor  \right)  }{f(i)} \right)^{ \frac{i}{\kappa}    } ( f(i))^y  =0. \label{eq:gen_condition_for_existance}
\end{equation}
Then, there exists a discrete probability distribution $P_W$ supported on the set $\cW $ defined in \eqref{eq:arb_support} such that
\begin{equation}
    \bbE[W^y 1_{W \le f(y)}] = \bbE[W^y 1_{W > f(y)}], \, y  \in \bbN_0. 
    \end{equation}
\end{theorem}

Our approach is to first show a truncated version, where the first $M$ moment conditions are satisfied, and then show that we can pass to a subsequence whose limit exists as the truncation point $M$ is moved to infinity. 

\subsection{Proof of  the Truncated Version of Theorem~\ref{thm:moment_theorem}}

We first consider a truncated version of our problem. That is, given $M \ge 1$, we seek to find $P_W$ such that 
\begin{equation}
    \bbE[W^y 1_{W \le f(y)}] = \bbE[W^y 1_{W > f(y)}], \, y  \in [0:M-1] .\label{eq:equivalent_cond_truncated}
    \end{equation}
 Naturally, we assume that $M\le c_0$. 

We start with the following helper lemma. 

\begin{lemma} \label{lem:helperlemma}
    For every $M \ge 1$ and $0< w_1 < \ldots <w_{M+1} <\infty$, there exists a probability vector $(p_1, \ldots, p_{M+1})$ such that  for all $k \in [0:M-1]$
    \begin{equation}
    \sum_{i=1}^{k+1} p_i w_i^k = \sum_{i=k+2}^{M+1} p_i w_i^k .
    \end{equation}
\end{lemma}
\begin{IEEEproof}
    We show this by induction.  For the base case of $M=1$, we have $ p_1  = p_2$,
    which results in valid probability vector if  $p_1 =  p_2  =1/2$. 

We now make an induction hypothesis that the statement is true for  $M$.    Let  $w_{M+2} > w_{M+1}$.  Now by the induction hypothesis for $(w_1,\ldots, w_{M+1})$ there exists a probability vector $(p_1, \ldots, p_{M+1})$ such that for all $k \in [0:M-1]$
\begin{equation}
\sum_{i=1}^{k+1} p_i w_i^k = \sum_{i=k+2}^{M+1} p_i w_i^k . 
\end{equation}
Moreover, by the induction hypothesis for $(w_1, \ldots, w_{M}, w_{M+2})$  there exists a probability vector $(q_1,\ldots, q_{M+1})$ such that  for all $k \in [0:M-1]$
\begin{equation}
\sum_{i=1}^{k+1} q_i w_i^k = \sum_{i=k+2}^{M} q_i w_i^k  + q_{M+1} w_{M+2}^k.\label{eq:equalities_for_qs}
\end{equation}

We now construct a solution to the case $M+1$ for the sequence $(w_1, \ldots, w_M, w_{M+1},w_{M+2})$. Fix some $\alpha  \in [0,1]$ (to be determined later) and let
\begin{align}
(m_1, \ldots, m_{M+2}) & = \alpha (p_1, \ldots, p_{M+1}, 0) \notag\\
&+ (1-\alpha) (q_1,\ldots,q_M, 0, q_{M+1} )
\end{align}
which obviously is a probability vector.  Moreover, by construction, for all $k \in [0:M-1]$, we have 
\begin{equation}
\sum_{i=1}^{k+1} m_i w_i^k = \sum_{i=k+2}^{M+2} m_i w_i^k  .
\end{equation}
To conclude the proof, we need to examine the case of  $k = M$ and show that
\begin{equation}
\sum_{i=1}^{M+1} m_i w_i^M =m_{M+2} w_{M+2}^M , 
\end{equation}
which can be re-written as 
\begin{equation}
   \alpha \sum_{i=1}^{M+1} p_i w_i^M + (1-\alpha)  \sum_{i=1}^{M} q_i w_i^M = (1-\alpha)  q_{M+1} w_{M+2}^M . \label{eq:Final_equation_of_lemma}
\end{equation}
To show that there exists a choice of $\alpha \in [0,1]$ such that \eqref{eq:Final_equation_of_lemma} is true note that for $\alpha=1$ the LHS of \eqref{eq:Final_equation_of_lemma} is larger than the RHS.  For $\alpha =0$ the LHS is smaller than RHS since 
\begin{equation}
\sum_{i=1}^{M} q_i w_i^M  \le w_{M+2 }  \sum_{i=1}^{M} q_i w_i^{M-1} =  w_{M+2 }  q_{M+1} w_{M+2}^{M-1}
\end{equation}
where the last equality follows from \eqref{eq:equalities_for_qs}. This  shows that there  the desired $\alpha$  exists and concludes the proof. 
\end{IEEEproof}

We now show that for every increasing function there  exists a $P_W$ that satisfies \eqref{eq:equivalent_cond_truncated}.    Our construction works by fixing $P_W$ to be a discrete distribution with the support given by
\begin{equation}
\cW^M = \left\{w_i:  w_i = f(i-1), \, i \in [1:M+1]  \right\}. \label{eq:w_i's}
\end{equation}
\begin{theorem} \label{thm:truncated_version} Suppose  that  $f: \bbN_0 \to \bbR_{+}$ is increasing up to $c_0$ and $1 \le M \le c_0$.  Then, there exists a probability distribution $P_W$ supported on $\cW^M$ that satisfies \eqref{eq:equivalent_cond_truncated}. 
\end{theorem}
\begin{IEEEproof}
    With the choice of $w_i$'s in \eqref{eq:w_i's}, the system in   \eqref{eq:equivalent_cond_truncated} can be written as: for all $y  \in [0:M-1] $  
\begin{equation}
\sum_{i=1}^{M+1} p_i w_i^y 1_{ \{w_i \le f(y) \}} = \sum_{i=1}^{M+1} p_i w_i^y 1_{ \{w_i > f(y) \}} ,
\end{equation}
and which can be further simplified to 
\begin{equation}
\sum_{i=1}^{y+1} p_i w_i^y  = \sum_{i= y+2}^{M+1} p_i w_i^y  , \, y  \in [0:M-1] .\label{eq:system_equation_p} 
\end{equation}

Now by invoking Lemma~\ref{lem:helperlemma}, we see that there exists  a probability vector $(p_1,\ldots,p_{M+1})$ such that the system of equalities in \eqref{eq:system_equation_p} is satisfied. 
This concludes the proof. 
\end{IEEEproof}

We now produce a concentration bound that will be useful.
\begin{lemma}\label{lem:tail_bound} Let $W$ be the random variable constructed in Theorem~\ref{thm:truncated_version}. Then, for $\kappa \ge 1$  and   $  \frac{i}{ \kappa } \le M +1 $, 
\begin{equation}
    P[W \ge w_{i+1}] \le \left(  \frac{f \left( \lfloor \frac{i}{\kappa}  \rfloor  \right)  }{f(i)} \right)^{ \frac{i}{\kappa}    }  . \label{eq:concentration_bound}
    \end{equation}
\end{lemma}
\begin{IEEEproof}
\begin{align}
P[ W \ge w_{ i+1} ]  
&= \bbE [ 1_{ \{ W \ge f(i) \} } ] \\
& \le \bbE \left[ \left( \frac{W}{f(i)} \right)^{ \lfloor \frac{i}{\kappa}  \rfloor  } 1_{ \{ W \ge f(i) \} } \right] \\
& \le \bbE \left[ \left( \frac{W}{f(i)} \right)^{ \lfloor \frac{i}{\kappa}  \rfloor  } 1_{ \{ W \ge  f \left( \lfloor \frac{i}{\kappa}  \rfloor  \right)  \} } \right] \label{eq:using_monotonicity}  \\
& =  \bbE \left[ \left( \frac{W}{f(i)} \right)^{ \lfloor \frac{i}{\kappa}  \rfloor  }  1_{ \{ W \le  f \left( \lfloor \frac{i}{\kappa}  \rfloor  \right)  \} } \right] \label{eq:using_moment_cond} \\
& \le  \left( \frac{f \left( \lfloor \frac{i}{\kappa}  \rfloor  \right)  }{f(i)} \right)^{ \lfloor \frac{i}{\kappa}  \rfloor  } ,
\end{align}    
where \eqref{eq:using_monotonicity} follows from the assumption that $f$ is non-decreasing; and  the equality in \eqref{eq:using_moment_cond} follows by using the fact that $P_W $ satisfies \eqref{eq:equivalent_cond_truncated}. 
\end{IEEEproof}

\subsection{Extension to the General Case}

 If $c_0<\infty$, then the proof is done; so  here, we consider the case of $c_0=\infty$. 
 For $M  \ge 1$, let $P_W^M$ denote a distribution constructed in Theorem~\ref{thm:truncated_version} for some $M>1$. In this section, we show that there is at least a subsequence $M_k$ such that $P_W = \lim_{k \to \infty} P_{W}^{M_k}$ converges (in total variation) and $P_W$ is a valid  distribution that  satisfies the desired system.

Let $p^M = (p^M_1,\ldots,p^M_{M+1})$ denote the probability vector corresponding to $P_W^M$. Using Lemma~\ref{lem:tail_bound} with $\kappa>1$, given $\varepsilon>0$, let $M_{\varepsilon,1},i$ be such that: for all $M\geq M_{\varepsilon,1}$
\begin{equation}P[W_M \geq  f(i)] \leq \frac{\varepsilon}{3}. \label{eq:epsilon_step1} \end{equation}

Let $\Pi_i(p^M)$ denote the projection onto the first $i$ components of $p^M$ (or equivalently setting all components $p_k^M$ with $k>i$ to zero). Using sequential compactness, let $M_k$ be a subsequence such that $\Pi_i(p^{M_k})$ converges in $\ell^1$ norm as $k\to\infty$. Since this subsequence is also Cauchy, let $M_{\varepsilon,2}$ be such that
\begin{equation} \|\Pi_i(p^{M_k}) - \Pi_i(p^{M_j}) \|_1 \leq \frac{\varepsilon}{3} \label{eq:epsilon_step2}\end{equation}
for all $M_j,M_k>M_{\varepsilon,2}$. Combining  \eqref{eq:epsilon_step1} and \eqref{eq:epsilon_step2}, we have
\begin{equation} \|p^{M_k} - p^{M_j}\|_1 \leq \varepsilon \end{equation}
for all $M_k,M_j \geq \max(M_{\varepsilon,1},M_{\varepsilon,2})$.

So far we have only constructed a subsequence that works for one particular $\varepsilon$, not any $\varepsilon>0$. This can be remedied by setting $\varepsilon_\ell = \frac{1}{2^\ell}$ and considering the following procedure. First, we construct a subsequence $M_k$ that works for $\varepsilon_1$. We leave this subsequence alone for $M_k \leq \max(M_{\varepsilon_1,1},M_{\varepsilon_1,2})$, but for $M_k > \max(M_{\varepsilon_1,1},M_{\varepsilon_1,2})$ we refine this subsequence by taking another subsequence of it that works for $\varepsilon_2$. Repeating this process results in a final Cauchy subsequence $M_k$, which by completeness also converges in $\ell^1$.

We now show that $P_W$ is a valid distribution.  Note that 
\begin{align}
\sum_{i=1}^\infty P_{W}(w_i)  &= \lim_{N \to \infty} \sum_{i=1}^N  P_{W}(w_i)\\
 &= \lim_{N \to \infty}  \lim_{k \to \infty} \sum_{i=1}^N  P_{W}^{M_k}(w_i) \le 1 .\label{eq:upper bound sum }
\end{align}
We now show that  $\sum_{i=1}^\infty P_{W}(w_k) \ge 1$.  First, note that for large enough $M$, from Lemma~\ref{lem:tail_bound}, we have that for all $k$
\begin{equation}
P^M[ W \ge w_{k+1}]   \le  \left(  \frac{f \left( \lfloor \frac{i}{\kappa}  \rfloor  \right)  }{f(i)} \right)^{ \frac{i}{\kappa}    } . \label{eq:concentration_bound}
\end{equation}

Now note that for all $j$
\begin{align}
\sum_{i=1}^{j+1}  P_{W}(w_i)& = \lim_{k \to \infty} \sum_{i=1}^{j+1}  P_{W}^{M_k}(w_i)\\
&=  \lim_{k \to \infty} (1- P[W_{M_k} \ge  w_{j+1}  ]  ) \\
& \ge 1- \left(  \frac{f \left( \lfloor \frac{j}{\kappa}  \rfloor  \right)  }{f(j)} \right)^{ \frac{j}{\kappa}    } , \label{eq:lower_bound_sum}
\end{align}
 where from \eqref{eq:gen_condition_for_existance} we have that 
$ \lim_{j \to \infty} \left(  \frac{f \left( \lfloor \frac{j}{\kappa}  \rfloor  \right)  }{f(j)} \right)^{ \frac{j}{\kappa}    } = 0$.

Combining \eqref{eq:upper bound sum } and \eqref{eq:lower_bound_sum}, and using  that $f$ is increasing,  we conclude that $P_{W}$ is a valid probability distribution. 

It remains to show that $P_W$ satisfies the following: for  $y \ge 0$
\begin{equation}
\sum_{i=1}^{y+1} P_{W}(w_i) w_i^y  = \sum_{i= y+2}^{\infty} P_{W}(w_i) w_i^y .
\end{equation}

Note that 
\begin{align}
    &\sum_{i= y+2}^{\infty} P_{W}(w_i) w_i^y  \notag\\
    &= \lim_{N \to \infty} \sum_{i= y+2}^{N} P_{W}(w_i) w_i^y  \\ 
    &= \lim_{N \to \infty}  \lim_{k \to \infty}  \sum_{i= y+2}^{N} P_{W}^{M_k}(w_i) w_i^y  \\ 
    &=   \lim_{k \to \infty}  \sum_{i= y+2}^{\infty} P_{W}^{M_k}(w_i) w_i^y  - \lim_{N \to \infty}  \lim_{k \to \infty}  \sum_{i= N+1}^{\infty} P_{W}^{M_k}(w_i) w_i^y   \\ 
    &=   \lim_{k \to \infty}  \sum_{i= y+2}^{\infty} P_{W}^{M_k}(w_i) w_i^y  - \Delta   \\ 
    &=   \lim_{k \to \infty}  \sum_{i= 1}^{y+1} P_{W}^{M_k}(w_i) w_i^y  - \Delta   \\ 
    &=    \sum_{i= 1}^{y+1} P_{W}(w_i) w_i^y  - \Delta    . 
\end{align}
Finally, note that
\begin{align}
    \Delta  &= \lim_{N \to \infty}  \lim_{k \to \infty}  \sum_{i= N+1}^{\infty} P_{W}^{M_k}(w_i) w_i^y \\
    &\le  \lim_{N \to \infty}  \sum_{i= N+1}^{\infty}   \left(  \frac{f \left( \lfloor \frac{i}{\kappa}  \rfloor  \right)  }{f(i)} \right)^{ \frac{i}{\kappa}    } ( f(i))^y  \\
    &=0, \label{eq:bounds_on_sum}
\end{align}
where \eqref{eq:bounds_on_sum} follows from the bound in \eqref{eq:concentration_bound} and the fact that the series is convergent from the assumption in \eqref{eq:gen_condition_for_existance}.

Therefore, we have that for all $y \ge  0$
\begin{equation}
    \sum_{i= y+2}^{\infty} P_{W}(w_i) w_i^y   
    =    \sum_{i= 1}^{y+1} P_{W}(w_i) w_i^y   . 
\end{equation}
This concludes the proof.

\subsection{Proof of Theorem~\ref{thm:main_theorem}}
Combining Lemma~\ref{lem:Moment_Version} and Theorem~\ref{thm:moment_theorem}, it remains to show that the $P_X$ is a valid probability distribution: 
    \begin{align}
    P_X(w_{i+1}) &\propto  P_W(w_{i+1}) \rme^{w_{i+1}} \\ 
    &\le  \left(  \frac{f \left( \lfloor \frac{i}{\kappa}  \rfloor  \right)  }{f(i)} \right)^{ \frac{i}{\kappa}    } \rme^{f(i)} ,\label{eq:used_tail_bound}
    \end{align}
    where \eqref{eq:used_tail_bound} follows from the tail bound in Lemma~\ref{lem:tail_bound}, which results in a valid distribution provided that the sum converges.
 The proof is concluded by noting that in Theorem~\ref{thm:main_theorem} we assume the series in \eqref{eq:used_tail_bound} is summable.  \hfill $\blacksquare$

 \section{Examples}
\label{sec:Examples}
The distribution in Theorem~\ref{thm:main_theorem} can be approximated to any degree of accuracy  by the following procedure:
\begin{enumerate}
\item Choose some $M \ge 1$ and define 
\begin{equation}
\bfA_M =\begin{pmatrix}
  1 & 1 & 1 & \cdots & 1 \\
-1 & 1 & 1 & \cdots & 1 \\
-w_1 & -w_2 & w_3 & \cdots & w_{M+1} \\
-w_1^2 & -w_2^2 & -w_3^2 & \cdots & w_{M+1}^2 \\
\vdots & \vdots & \vdots & \ddots & \vdots \\
-w_1^{M-1} & -w_2^{M-1} & -w_3^{M-1} & \cdots & w_{M+1}^{M-1}
\end{pmatrix} 
\end{equation}
where $w_i = f(i-1),  \, i \in [1:M+1]$ and where $f$ is the desired estimator. 
\item Let
\begin{equation}
\bfp_M =  \bfA_M^{-1} \begin{pmatrix}
1 , 0, \ldots, 0
\end{pmatrix} ^T
\end{equation}
where $\bfA_M$ is guaranteed to have an inverse. 
\item A truncated approximation of $P_X$ is given by: for $i \in [1:M+1]$
\begin{equation}
    P_X^{(M)}(w_i) \propto \bfp_M(i) \rme^{w_i}.   \label{eq:construcion_alg}
\end{equation}
\end{enumerate}

\begin{figure}
    \centering  \input{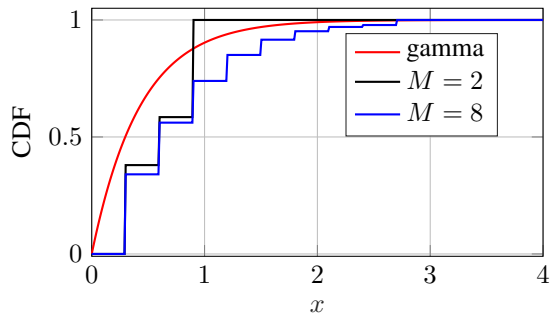}
    \caption{Comparison of the gamma cdf (red),  which induces a linear conditional mean, to the distribution in \eqref{eq:construcion_alg}. }
    \label{fig:cdfs}
\end{figure}

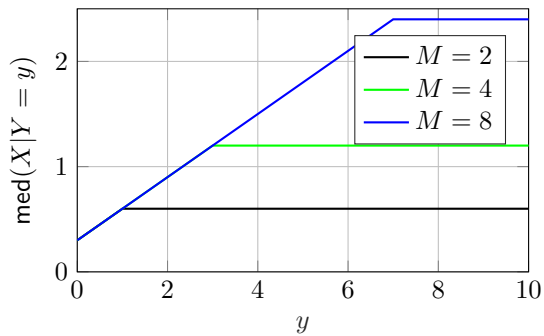
\begin{figure}
    \centering  
%
%
\definecolor{mycolor1}{rgb}{0.00000,0.44700,0.74100}%
\definecolor{mycolor2}{rgb}{0.85000,0.32500,0.09800}%
\definecolor{mycolor3}{rgb}{0.92900,0.69400,0.12500}%
\begin{tikzpicture}

\begin{axis}[%
width=6cm,
height=3.5cm,
at={(1.011in,0.642in)},
scale only axis,
xmin=0,
xmax=10,
xlabel style={font=\color{white!15!black}},
xlabel={$y$},
ylabel={${\sf med}(X|Y=y)$},
ymin=0,
ymax=2.5,
axis background/.style={fill=white},
xmajorgrids,
ymajorgrids,
legend style={legend cell align=left, align=left, draw=white!15!black, at = {(0.95,0.9)}}
]
\addplot [color=black, thick]
  table[row sep=crcr]{%
0	0.3\\
1	0.6\\
2	0.6\\
3	0.6\\
4	0.6\\
5	0.6\\
6	0.6\\
7	0.6\\
8	0.6\\
9	0.6\\
10	0.6\\
};
\addlegendentry{$M=2$}

\addplot [color=green, thick]
  table[row sep=crcr]{%
0	0.3\\
1	0.6\\
2	0.9\\
3	1.2\\
4	1.2\\
5	1.2\\
6	1.2\\
7	1.2\\
8	1.2\\
9	1.2\\
10	1.2\\
};
\addlegendentry{$M=4$}

\addplot [color=blue, thick]
  table[row sep=crcr]{%
0	0.3\\
1	0.6\\
2	0.9\\
3	1.2\\
4	1.5\\
5	1.8\\
6	2.1\\
7	2.4\\
8	2.4\\
9	2.4\\
10	2.4\\
};
\addlegendentry{$M=8$}

\end{axis}

\end{tikzpicture}%
    \caption{Conditional medians induced by the distribution in \eqref{eq:construcion_alg} with $a=b=0.3$.}
    \label{fig:cond_medians}
\end{figure}

It is interesting to simulate the truncated approximation in \eqref{eq:construcion_alg} for the affine case:   

\begin{itemize}
\item Fig.~\ref{fig:cdfs} shows the cdf of the distribution in \eqref{eq:construcion_alg} for $f(y) =a y +b$ with $a=b=0.3$ and  several values of $M$ and compares it to the cdf of the gamma distribution. 
\item Fig.~\ref{fig:cond_medians} shows conditional medians of the distribution in \eqref{eq:construcion_alg} for $M=2, 4$, and $8$ with $a=b=0.3.$
\end{itemize}

\section{Conclusion}
This work has addressed the problem of $L^1$ estimation under Poisson noise, focusing on the Bayesian case. By constructing a novel family of priors, we have demonstrated that any increasing function can serve as an admissible Bayesian estimator, provided certain integrability conditions are met.

In particular, we highlighted a class of priors distinct from the classical gamma distribution, which is known to induce affine conditional means under $L^2$ loss. These new priors represent a significant departure from the conjugate priors typically used for exponential families and are constructed using a limiting process that guarantees the desired linearity properties.

Future work could explore extensions of this framework to other noise models such as Binomial or Negative Binomial distributions, or to different loss functions. 
\bibliography{refs.bib}
 \bibliographystyle{IEEEtran}
\end{document}